\documentclass [11pt]{article}
\usepackage{amsthm}
\usepackage{amsmath}
\usepackage{amssymb}
\usepackage{graphics}
\usepackage{epsfig}
\usepackage{color}

\newcounter{contador}
\newtheorem{propo}[contador]{Proposition}
\newtheorem{teo}[contador]{Theorem}

\newtheorem{nota}[contador]{Remark}

\newcommand{\rec}{\noindent}    

\newcommand{\g}{\gamma}
\newcommand{\G}{\Gamma}
\newcommand{\tG}{ \widetilde{\Gamma}}

\newcommand{\enya}{${\rm \tilde{n}}$}

\newcommand{\tp}{ \widetilde{P}}
\newcommand{\ag}{\operatorname{Arg}}

\newcommand{\R}{{\mathbb R}}
\newcommand{\C}{{\mathbb C}}

\newcommand{\U}{{\cal{U}}}

\title{On  Poncelet's Maps}

\author{Anna Cima$^{(1)}$, Armengol Gasull$^{(1)}$ and V\'{\i}ctor Ma\~{n}osa $^{(2)}$
  \\*[.1truecm]
{\small \textsl{$^{(1)}$ Dept. de Matem\`{a}tiques, Facultat de
Ci\`{e}ncies,}}
\\*[-.25truecm] {\small \textsl{Universitat Aut\`{o}noma de Barcelona,}}
\\*[-.25truecm] {\small \textsl{08193 Bellaterra, Barcelona, Spain}}
\\*[-.25truecm] {\small \textsl{cima@mat.uab.cat, gasull@mat.uab.cat}}
\\*[-.25truecm]
\\*[-.25truecm] {\small \textsl{$^{(2)}$ Dept. de Matem\`{a}tica Aplicada III (MA3),}}
\\*[-.25truecm] {\small \textsl{Control, Dynamics and Applications Group (CoDALab)}}
\\*[-.25truecm] {\small \textsl{Universitat Polit\`{e}cnica de Catalunya (UPC)}}
\\*[-.25truecm] {\small \textsl{Colom 1, 08222 Terrassa, Spain}}
\\*[-.25truecm] {\small \textsl{victor.manosa@upc.edu}}}
\date{}

\begin{document}
\maketitle

\begin{abstract}
Given two ellipses, one surrounding the other one, Poncelet
introduced a map $P$ from the exterior one to itself by using the
tangent lines to the interior ellipse. This procedure can be
extended to any two smooth, nested and  convex  ovals and we call
this type of maps Poncelet's maps. We recall what he proved around
1814 in the dynamical systems language: In the two ellipses case and
when the rotation number of $P$  is rational there exists a
$n\in\mathbb{N}$ such that $P^n=\operatorname{Id},$ or in other
words,   the Poncelet's map is conjugated to a rational rotation. In
this paper we study general Poncelet's maps and give several
examples of algebraic ovals where the corresponding Poncelet's map
has a rational rotation number and it is not conjugated to a
rotation. Finally, we also provide a new proof of Poncelet's result
based on dynamical tools.
\end{abstract}

\rec {\sl 2000 Mathematics Subject Classification:} Primary: 51M15.
Secondary: 37C05, 51M04.

 \rec {\sl Keywords:} Poncelet's problem, circle maps, rotation number, devil's staircase\newline


\section{Introduction and Main Results}\label{sec1}

Let $\gamma$ and $\Gamma$ be two $\mathcal{C}^r,$ $r\ge1,$ simple,
closed and nested curves, each one of them  being the boundary of  a
convex set. Furthermore we assume for instance that $\Gamma$
surrounds~$\gamma.$

Given any $p\in \Gamma$ there are exactly two points $q_1, q_2$ in
$\gamma$ such that the lines $p\,q_1,$ $p\,q_2$ are tangent to
$\gamma.$ We define the {\it Poncelet's map},
$P:\Gamma\rightarrow\Gamma,$ associated to the pair $\gamma,\Gamma$
as
\[
P(p)=P_{\gamma,\Gamma}(p)=\overline{pq_1}\cap \Gamma,
\]
where $p\in\G,$ $\overline{pq_1}\cap \Gamma$ is the first point in
the set $\{\overline{pq_1}\cap \Gamma, \overline{pq_2}\cap \Gamma\}$
that we find when, starting from $p,$ we follow $\Gamma$
counterclockwise, see  Figure 1. Notice that $
P^{-1}(p)=\overline{pq_2}\cap \Gamma$.


\begin{center}
\begin{picture}(0,0)%
\epsfig{file=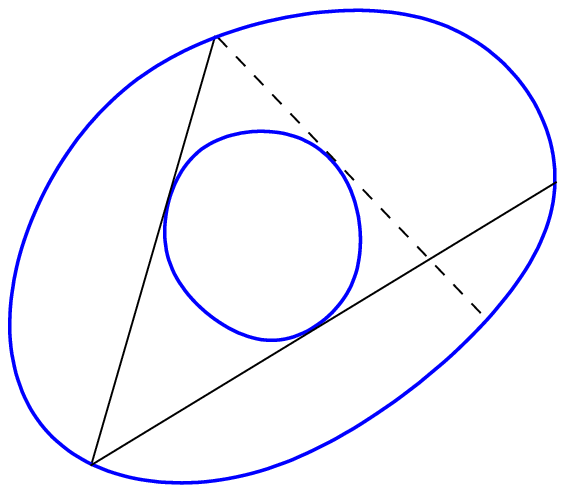}
\end{picture}%
\setlength{\unitlength}{4144sp}%
\begingroup\makeatletter\ifx\SetFigFont\undefined%
\gdef\SetFigFont#1#2#3#4#5{%
  \reset@font\fontsize{#1}{#2pt}%
  \fontfamily{#3}\fontseries{#4}\fontshape{#5}%
  \selectfont}%
\fi\endgroup%
\begin{picture}(3350,2487)(478,-1918)
\put(1519,395){\makebox(0,0)[rb]{\smash{{\SetFigFont{10}{12.0}{\rmdefault}{\mddefault}{\updefault}{\color[rgb]{0,0,0}$p$}%
}}}}
\put(925,-1873){\makebox(0,0)[rb]{\smash{{\SetFigFont{10}{12.0}{\rmdefault}{\mddefault}{\updefault}{\color[rgb]{0,0,0}$P(p)$}%
}}}}
\put(2869,-1177){\makebox(0,0)[lb]{\smash{{\SetFigFont{10}{12.0}{\rmdefault}{\mddefault}{\updefault}{\color[rgb]{0,0,0}$P^{-1}(p)$}%
}}}}
\put(3229,-373){\makebox(0,0)[lb]{\smash{{\SetFigFont{10}{12.0}{\rmdefault}{\mddefault}{\updefault}{\color[rgb]{0,0,0}$P^2(p)$}%
}}}}
\put(2719,461){\makebox(0,0)[lb]{\smash{{\SetFigFont{10}{12.0}{\rmdefault}{\mddefault}{\updefault}{\color[rgb]{0,0,0}$\Gamma$}%
}}}}
\put(1561,-907){\makebox(0,0)[lb]{\smash{{\SetFigFont{10}{12.0}{\rmdefault}{\mddefault}{\updefault}{\color[rgb]{0,0,0}$\gamma$}%
}}}}
\put(1243,-517){\makebox(0,0)[rb]{\smash{{\SetFigFont{10}{12.0}{\rmdefault}{\mddefault}{\updefault}{\color[rgb]{0,0,0}$q_1$}%
}}}}
\put(2173,-241){\makebox(0,0)[lb]{\smash{{\SetFigFont{10}{12.0}{\rmdefault}{\mddefault}{\updefault}{\color[rgb]{0,0,0}$q_2$}%
}}}}
\end{picture}%
\\
\textbf{FIGURE 1.} The Poncelet's map.
\\
\end{center}


The implicit function theorem together with the geometrical
interpretation of the construction of $P$ imply that it is a
$\mathcal{C}^r$ diffeomorphism from $\Gamma$ into itself. So $P$ can
be seen as a $\mathcal{C}^r$ diffeomorphism of the circle and has
associated a {\it rotation number}
\[
\rho=\rho(P)=\rho(\gamma,\Gamma)\in(0,1/2).
\]
See for instance   \cite{ALM,AP} for the definition of rotation
number. Notice that usually a rotation number is in $(0,1).$ Our
choice of $q_1$ for the Poncelet's map implies that indeed
$\rho<1/2.$ It is also well known that if $\Phi$ is any
diffeomorphism of the circle of class at least $\mathcal{C}^2$ and
such that $\rho(\Phi)\not\in\mathbb{Q}$ then $\Phi$ is conjugated to
a rotation of angle $\rho(\Phi).$  So this is the situation for the
Poncelet's map $P$ when $\rho(P)\not\in\mathbb{Q}$ and $r\ge2.$

With the above notation the celebrated {\it Poncelet's Theorem}
asserts that if $\gamma$ and $\Gamma$ are ellipses, with arbitrary
relative positions,  and $\rho=\rho(\gamma,\Gamma)\in \mathbb{Q}$
then the Poncelet's map is also conjugated to the rotation of angle
$\rho$ in $\mathbb{S}^1.$ In geometrical terms, if  starting at some
point $p\in\Gamma$ the Poncelet's process of drawing tangent lines
to $\gamma$ closes after $n$ steps then the same holds for any other
starting point in $\Gamma.$ There are several proofs of this nice
result  in \cite[Sec. 4.3]{T} and a different one, based on a
beautiful approach of Bertrand and Jacobi through differential
equations and elliptic integrals in \cite[pp. 191-194]{S}. In
Section \ref{newproof} we give another proof based on dynamical
tools, by using the results of \cite{CGM}.  The problem of
determining explicit conditions over the coefficients of the two
ellipses to ensure that the Poncelet's map is conjugated to a
rational rotation was solved by Cayley. An excellent exposition of
this result is given in~\cite{GH}.

A  monograph devoted to  Poncelet's theorem and related results it
is going to appear soon, see \cite{F}.

It is clear that in Poncelet's Theorem it is not restrictive to
assume that $\gamma=\{x^2+y^2=1\}.$ The first question that we face
in this paper is the following:  Is the Poncelet's result also true
when we consider $\gamma=\{x^2+y^2=1\}$ and $\Gamma$ given by an
oval of an algebraic curve of higher degree? We prove:

\begin{teo}\label{mt} Fix
$\gamma=\{x^2+y^2=1\}.$ Then for any $m\in \mathbb{N,}$  $m>2,$
there is an algebraic curve of degree $m,$  containing  a convex
oval $\G,$ such that   the Poncelet's map associated to $\g$ and
$\G$ has rational rotation number and it is not conjugated to  a
rotation.
\end{teo}

This result will be  a consequence of a more general result proved
in Section \ref{sec2}, see Proposition \ref{prop3}. Moreover the
full dynamics of the introduced Poncelet's maps $P:\Gamma\rightarrow
\Gamma$ will be described in that section.

From  Theorem \ref{mt} it is clear that, in general,  Poncelet's
maps with rational rotation numbers are not  conjugated to
rotations. It is natural to  wonder about this question  when both
ovals are level sets of the same polynomial map
$V:\mathbb{R}^2\rightarrow \mathbb{R}.$ As one of the simplest cases
we consider the homogeneous map $V(x,y)=x^4+y^4$, giving
\[
\gamma=\Gamma_1=\{x^4+y^4=1\} \quad \mbox{and}  \quad
\Gamma=\Gamma_k=\{x^4+y^4=k\},
\]
for $k\in\mathbb{R}, k>1.$ As we will see in  Section \ref{sec3}
even in this  situation the conjugacy with a rotation is not
 true.

The Poncelet's maps also provide  a natural way of defining an
integrable function from an open set of $\mathbb{R}^2$ into itself
as follows: We foliate the open unbounded connected component
$\mathcal{U}$ of $\mathbb{R}^2\setminus \Gamma_1$ as
\begin{equation*}\label{open}
\mathcal{U}:=\bigcup_{k>1}\{x^4+y^4=k\}.
\end{equation*}
Then  the Poncelet's construction gives a new diffeomorphism, also
of class $\mathcal{C}^r,$ that is defined from $\mathcal{U}\subset
\mathbb{R}^2$ into itself, which simply consists in applying to each
point $p$ the corresponding Poncelet's map, associated to the level
set of $V$ passing trough $p.$ For sake of simplicity we also call
it~$P.$ Notice that this map is trivially integrable by means of
$V(x,y)=x^4+y^4,$ that is $V(P(x,y))=V(x,y)$ for all
$(x,y)\in\mathcal{U}.$

As we will see in Subsection \ref{ode}, this extended Poncelet's map
$P$ will  be useful to give properties of a  functional equation
that helps to study integrable planar maps.

Finally, in the above context it is natural to introduce the
rotation function $\rho(k):=\rho(\Gamma_1,\Gamma_k),$ as the
rotation number of $P$ associated to $\gamma$ and $\Gamma_k,$ and to
study some of its properties.

In the case of two  concentric circles
\[
\gamma=\widetilde\Gamma_1=\{x^2+y^2=1\} \quad \mbox{and}  \quad
\Gamma=\widetilde\Gamma_k=\{x^2+y^2=k\},
\]
it is easy to prove that the rotation function
$\widetilde\rho(k)=\rho(\widetilde\Gamma_1,\widetilde\Gamma_k)$ is
the smooth monotonous function $
\widetilde\rho(k)=\arctan(k-1)/{\pi}.
 $
On the other hand, in Section \ref{sec3} we will show that the
function $\rho(k):=\rho(\Gamma_1,\Gamma_k),$ is much more
complicated. Indeed all what we prove seems to indicate that it has
the usual shape of the rotation function of generic one-parameter
families of diffeomorphisms: the devil's staircase, see for instance
\cite{BG,O}.

\section{A Circle and an Oval}\label{sec2}

We prove a preliminary result that implies Theorem~\ref{mt}.

\begin{propo}\label{prop3} Consider
\[
\gamma=\{x^{2n}+y^{2n}=1\}  \quad \mbox{and}  \quad
\Gamma=\{x^{2m}+y^{2m}=2\}\quad \mbox{with}\quad n,m\in\mathbb{N},
\]
and let $P$ be the Poncelet's map associated to them. Then
$\rho_{n,m}(P)=1/4.$ Moreover, the map is conjugated a rotation if
and only if $n=m=1.$
\end{propo}

\begin{proof}
It is easy to check  that for any $n$ and $m$, the Poncelet map $P$
has the periodic orbit of period 4, given by the points
$\mathcal{O}=\{(1,1),(-1,1),(-1,-1),(1,-1)\}$. Hence
$\rho_{n,m}(P)=1/4.$

When $n=m=1,$  both ovals are coniques and Poncelet's  Theorem
proves one implication of our result. Let us prove the converse.
Assume that $P$ is conjugated to a rotation. Then, since
$\rho(P)=1/4,$ for each $p\in\G,$ $P^4(p)=p.$ Consider the
particular point $p=p_1:=(0,\sqrt[2m]{2})\in\G.$ By using the
symmetries of the problem, the 4-periodicity of $P$ forces that
$p_2:=P(p_1)=(-\sqrt[2m]{2},0)$. Thus the line trough $p_1$ and
$p_2$ has to be tangent to the oval $x^{2n}+y^{2n}=1$ at some point
$\bar p=(\bar x, \bar y).$ Putting all the conditions together gives
that $\bar p$  has to be solution of the system
$$
\begin{cases}\begin{array}{l}
         y=x+\sqrt[2m]{2}, \\
         x^{2n}+y^{2n}=1, \\
         (2nx^{2n-1},2ny^{2n-1})\cdot(1,1)=x^{2n-1}+y^{2n-1}=0.
       \end{array}
\end{cases}
$$
It has a real solution only when
\begin{equation}\label{con}
2n=\frac{2m}{2m-1}=1+\frac{1}{2m-1},
\end{equation}
and in this case it is  $\bar p=(-\sqrt[2m]{2}/2,\sqrt[2m]{2}/2).$
It is clear from \eqref{con} that the only  solution with natural
values of this equation is $n=m=1,$ as we wanted to
prove.\end{proof}

\begin{proof}[Proof of Theorem \ref{mt}] Clearly the proof when $m$
is even is a corollary of the Proposition~\ref{prop3}. The proof
when $m\ge3$ is odd follows by noticing that the sets
$\{x^{2m}+y^{2m}-2=0\}$  and
$\{(x+10)\left(x^{2m}+y^{2m}-2\right)=0\}$ coincide in $\{x>-10\}$
and so in both cases the Poncelet's maps coincide.
\end{proof}

For the simplest case studied in Proposition~\ref{prop3}, $n=1,$ we
prove the following result that characterizes the dynamics of $P.$

\begin{propo}\label{prop4} Consider
\[
\gamma=\{x^{2}+y^{2}=1\}  \quad \mbox{and}  \quad
\Gamma=\{x^{2m}+y^{2m}=2\}\quad \mbox{with}\quad m\in\mathbb{N}, m>
1.
\]
and let $P$ be the Poncelet's map associated to them. Then
$\rho(P)=1/4,$ the orbit $\mathcal{O}
=\{(1,1),(-1,1),(-1,-1),(1,-1)\}$ is a 4-periodic orbit of $P$ and
it is the $\alpha$ and the $\omega$ limit of all the orbits of $P.$
\end{propo}

\begin{proof}
Set  $\tG:=\{x^{2}+y^{2}=2\}$, and let $P$ and $\tp$ be the the
Poncelet maps associated to $\G$ and $\g$,  and  to $\tG$ and $\g,$
respectively. As usual  given a point $p_1\in \G$  (resp.
$q_1\in\tG$) we write $p_{i+1}:=P(p_i)\in\G$ (resp.
$q_{i+1}:=\tp(q_i)\in\tG$), for $i\ge1.$  Our goal will be to
compare both maps. In fact, we will use the map $\tp$ as a kind of
Lyapunov function for the map $P.$

\begin{center}
\begin{picture}(0,0)%
\epsfig{file=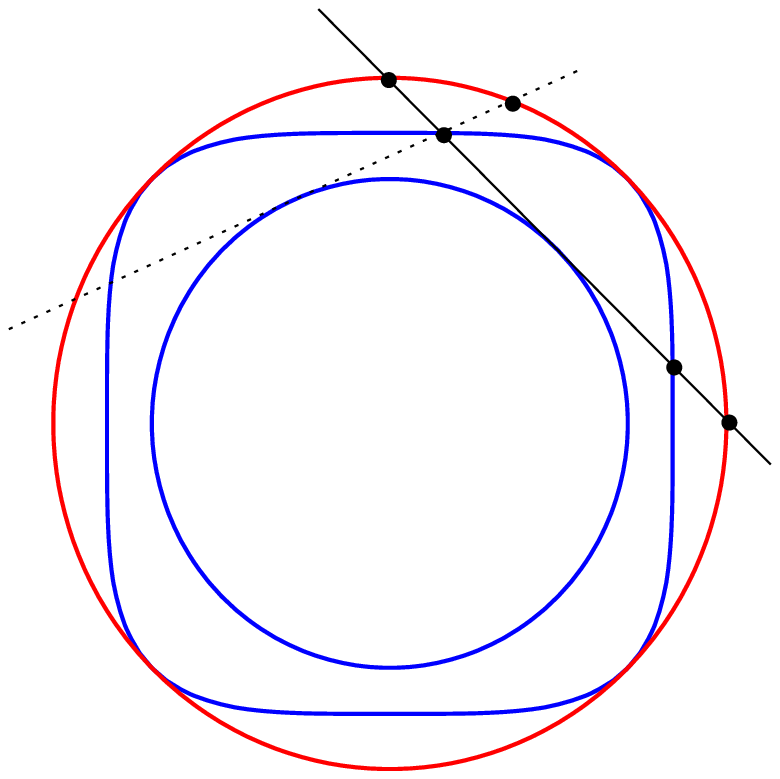}%
\end{picture}%
\setlength{\unitlength}{4144sp}%
\begingroup\makeatletter\ifx\SetFigFont\undefined%
\gdef\SetFigFont#1#2#3#4#5{%
  \reset@font\fontsize{#1}{#2pt}%
  \fontfamily{#3}\fontseries{#4}\fontshape{#5}%
  \selectfont}%
\fi\endgroup%
\begin{picture}(3600,3820)(451,-3431)
\put(2701,-376){\makebox(0,0)[lb]{\smash{{\SetFigFont{10}{12.0}{\rmdefault}{\mddefault}{\updefault}{\color[rgb]{0,0,0}$p_2$}%
}}}}
\put(3916,-1591){\makebox(0,0)[lb]{\smash{{\SetFigFont{10}{12.0}{\rmdefault}{\mddefault}{\updefault}{\color[rgb]{0,0,0}$q_1= \phi(p_1)$}%
}}}}
\put(3016,-2716){\makebox(0,0)[rb]{\smash{{\SetFigFont{10}{12.0}{\rmdefault}{\mddefault}{\updefault}{\color[rgb]{0,0,0}$\Gamma$}%
}}}}
\put(2071,-3391){\makebox(0,0)[lb]{\smash{{\SetFigFont{10}{12.0}{\rmdefault}{\mddefault}{\updefault}{\color[rgb]{0,0,0}$\widetilde \Gamma$}%
}}}}
\put(1576,-2266){\makebox(0,0)[lb]{\smash{{\SetFigFont{10}{12.0}{\rmdefault}{\mddefault}{\updefault}{\color[rgb]{0,0,0}$\gamma$}%
}}}}
\put(2026, 74){\makebox(0,0)[rb]{\smash{{\SetFigFont{10}{12.0}{\rmdefault}{\mddefault}{\updefault}{\color[rgb]{0,0,0}$q_2=\widetilde P (q_1)$}%
}}}}
\put(3016,-151){\makebox(0,0)[lb]{\smash{{\SetFigFont{10}{12.0}{\rmdefault}{\mddefault}{\updefault}{\color[rgb]{0,0,0}$\phi(p_2)$}%
}}}}
\put(3601,-1231){\makebox(0,0)[lb]{\smash{{\SetFigFont{10}{12.0}{\rmdefault}{\mddefault}{\updefault}{\color[rgb]{0,0,0}$p_1$}%
}}}}
\end{picture}%
\\
\textbf{FIGURE 2.} A comparison between the  Poncelet's maps $P$ and
$\tp$.
\end{center}

 Recall that Poncelet's Theorem  establishes that
$\tp$ is conjugated to a rotation. It is easy to check that the
orbit $\mathcal{O} =\{(1,1),(-1,1),(-1,-1),(1,-1)\}$ is a 4-periodic
orbit for $\tp$. So, $\rho(\tp)=1/4$ and then $\tp^4(q)=q$ for all
$q\in\tG$.

Notice that $\mathcal{O} $ is also a 4-periodic orbit of $P.$ Hence
$\rho(P)=1/4.$

Let us introduce some notation. Fixed $\g$ it is easy  to construct
a bijection between $\G$ and $\tG$ as follows: Given $p_1\in\G$ we
define $\phi(p_1)$ as the point of intersection between $\tG$  and
the half-straight line starting at $p_2=P(p_1)$ and passing trough
$p_1,$ see Figure 2. Notice that by construction $q_1:=\phi(p_1),$
$p_1$, $p_2$ and $q_2=\tp(q_1)$ are aligned. Given a point $r\ne0,$
we denote by $\ag(r)$ the argument modulus $2\pi$ of $r$ thought as
a point of $\C \setminus \{0\}$.

Take any $p_1\in\G\cap\{(x,y)\,:\, x>0, y>0$ and $y\neq x\}$, and
$q_1\in \tG$ given by $q_1=\phi(p_1)$. By construction $p_1$, $q_1$,
$p_2$ and $q_2$ are aligned. Notice also  that $ \ag(q_2) >
\ag(\phi(p_2)) $ where both angles are in $\left(0,\pi\right]$. This
can be understood as a ``delay'' of $P$ with respect to $\tp$. This
delay is propagated through the next three iterates giving that
$P^4$ does not complete a tour around $\G.$ Hence the lifting of the
map $P^4$ is below the identity except at the four points
corresponding to the 4-periodic orbit of $P.$ Then the result
follows.
\end{proof}

\begin{nota} From the proof of Proposition \ref{prop3} it is not
difficult to construct $\mathcal{C}^1$ ovals for which the dynamic
of the Poncelet's map is different from the one described in
Proposition~\ref{prop4}. By using the compatibility condition
\eqref{con}, it is easy to prove that, for instance when $n=2$ and
$m=2/3$ in the statement of Proposition \ref{prop3}, i.e. when
\[
\gamma=\{x^{4}+y^{4}=1\}  \quad \mbox{and}  \quad
\Gamma=\{x^{4/3}+y^{4/3}=2\},
\]
the Poncelet's map $P$ is such that $\rho(P)=1/4.$  It has  the two
4-periodic orbits $\mathcal{O}=\{(1,1),(-1,1),(-1,-1),(1,-1)\}$, and
$\mathcal{O}^*=\{(0,a),(-a,0),(0,-a),(a,0)\}$ where $a=2^{3/4}$ and
it is not conjugated to a rotation.\end{nota}

\section{Two Ovals}\label{sec3}

This section is devoted to study in more detail the Poncelet's maps
associated to the ovals $\gamma=\{x^4+y^4=1\}$ and
$\Gamma_k=\{x^4+y^4=k\},$ $k>1.$

\subsection{How to find periodic orbits}\label{po}

Let us impose which condition has to satisfy  an orbit to be $n$
periodic for a Poncelet's map associated to the two ovals
$\gamma=\{g(x,y)=0\}$ and $\Gamma=\{G(x,y)=0\}.$ Take $n$
counterclockwise  ordered points on $\gamma,$  $p_1,p_2,\ldots,p_n.$
Draw the $n$ tangent lines to $\gamma,$ corresponding to these
points, and denote them by $\ell_1,\ell_2,\ldots,\ell_n,$
respectively. To generate a $n$-periodic orbit of $P$ the following
assertions must hold:

\begin{itemize}
\item[(i)] The lines $\ell_i$ and $\ell_{i+1}$ for $i=1,2,\ldots, n,$ where $\ell_{n+1}:=\ell_{1},$
intersect. We denote these intersections by:
$q_{i,i+1}=\ell_i\cap\ell_{i+1}.$

\item[(ii)]\label{iitem} The following $n$ equalities are satisfied:
\[
G(q_{1,2})=G(q_{2,3})=\cdots=G(q_{n-1,n})=G(q_{n,n+1})=0.
\]

\end{itemize}

Notice that the above set of conditions gives a non-linear system
with $n$ unknowns (corresponding to the $n$  points on $\gamma$) and
$n$ equations.

On the other hand when we consider the same problem between
$\gamma=\{g(x,y)=0\}$ and $\Gamma_k=\{G(x,y)=k\}$ and we take $k$
also as a free unknown, the problem of searching $n$-periodic orbits
is equivalent to impose instead of item (ii), the following $n-1$
equalities:
\[
G(q_{1,2})=G(q_{2,3})=\cdots=G(q_{n-1,n})=G(q_{n,n+1}).
\]

So in this case we get again a  non-linear system with $n$ unknowns
 but now with only  $n-1$ equations. Hence,  it is natural to believe that either it
  has no solution or it has a continuum of them. Notice that this
  continuum can be interpreted as a continuum of values $k\in J_k\subset \mathbb{R}$  for
  which P has a periodic orbit of period $n.$ Observe also that each
  one of these continua gives rise to an interval $J_k$ where the rotation
  function $\rho(k)$ associated to the Poncelet's map between $\gamma=\{g(x,y)=0\}$ and
  $\Gamma_k=\{G(x,y)=k\}$ has a constant value $j/n,$ for some
  $j\in\mathbb{N}.$ Each  of these intervals will give rise to one
  stair in the devil's staircase which seems to be associated to $P.$

 In next subsections we  study in more detail the case
 $g(x,y)=x^4+y^4-1$ and $G(x,y)=x^4+y^4,$ and we will give a geometrical
 interpretation of the localization of the starting and the ending points of some of the  stairs of $\rho(k)$.

\subsection{Some symmetric periodic orbits}\label{spo}

 Notice that the ovals $\Gamma_k=\{x^4+y^4=k\}, k\ge 1$   have
several symmetries. These are given by the two axes and the
diagonals $\{(x,y)\in \mathbb{R}^2\,:\, y=\pm x\}.$ In this section
we search some periodic orbits of the Poncelet's maps $P,$
associated to $\Gamma_1$ and $\Gamma_k$ that share some of these
symmetries. For these type of periodic orbits the order of the
system described in the previous subsection, that has to be solved
to find the periodic orbits, can be reduced.

As an example we find a value $k$ for which the corresponding
Poncelet's map has rotation number $1/3$ due to the existence of a
symmetric 3-periodic orbit with respect to the $Oy-$axis. First we
take $\{(x_1,y_1),(0,1)\}\in\gamma.$ The corresponding tangent lines
to $\gamma$ are $\ell_1(x,y)=x_1^3x+y_1^3y-1=0$ and
$\ell_2(x,y)=y-1=0,$ see the left picture in Figure~3. We have that
$\ell_1\cap\ell_{2}=((1-y_1^3)/x_1^3,1).$ To give rise to a
3-periodic orbit with the searched symmetry, the third point has to
be $(x_3,y_3):=(-x_1,y_1).$ Moreover we get that
$\ell_1\cap\ell_{3}=\ell_1\cap\{x=0\}=(0,1/y_1^3).$ Hence the
conditions that imply that the three intersection points between the
tangent lines belong to the same  $\Gamma_k$  reduce to the single
equation
\[
\left(\frac{1-y_1^3}{x_1^3}\right)^4+1=
\left(\frac1{y_1^3}\right)^4,
\]
or equivalently to impose that $y_1$ satisfies the equation
$$R(y):=y^{12}(1-y^4)^3+y^{12}(1-y^3)^{4}-(1-y^4)^3=0,$$
where we have used that $x_1^4+y_1^4=1.$ Some calculations give that
\begin{eqnarray*}
R(y)&=&-4y^{21}+3y^{20}+6y^{18}-3y^{16}-4y^{15}+3y^{12}-3y^8+3y^4-1=\\&&-(y-1)^4(y^2+y+1)R_{15}(y),
\end{eqnarray*}
where $R_{15}$ is a polynomial of degree 15 which has only one real
root, $y=y_1^*\simeq -0.779644.$ The value of $k$ corresponding to
this solution is $k=k^*:=1/(y_1^*)^{12}\simeq19.8264.$ This result
is reflected in Table~1.

\begin{center}{
\begin{figure}[h]
$\qquad$ \epsfig{file=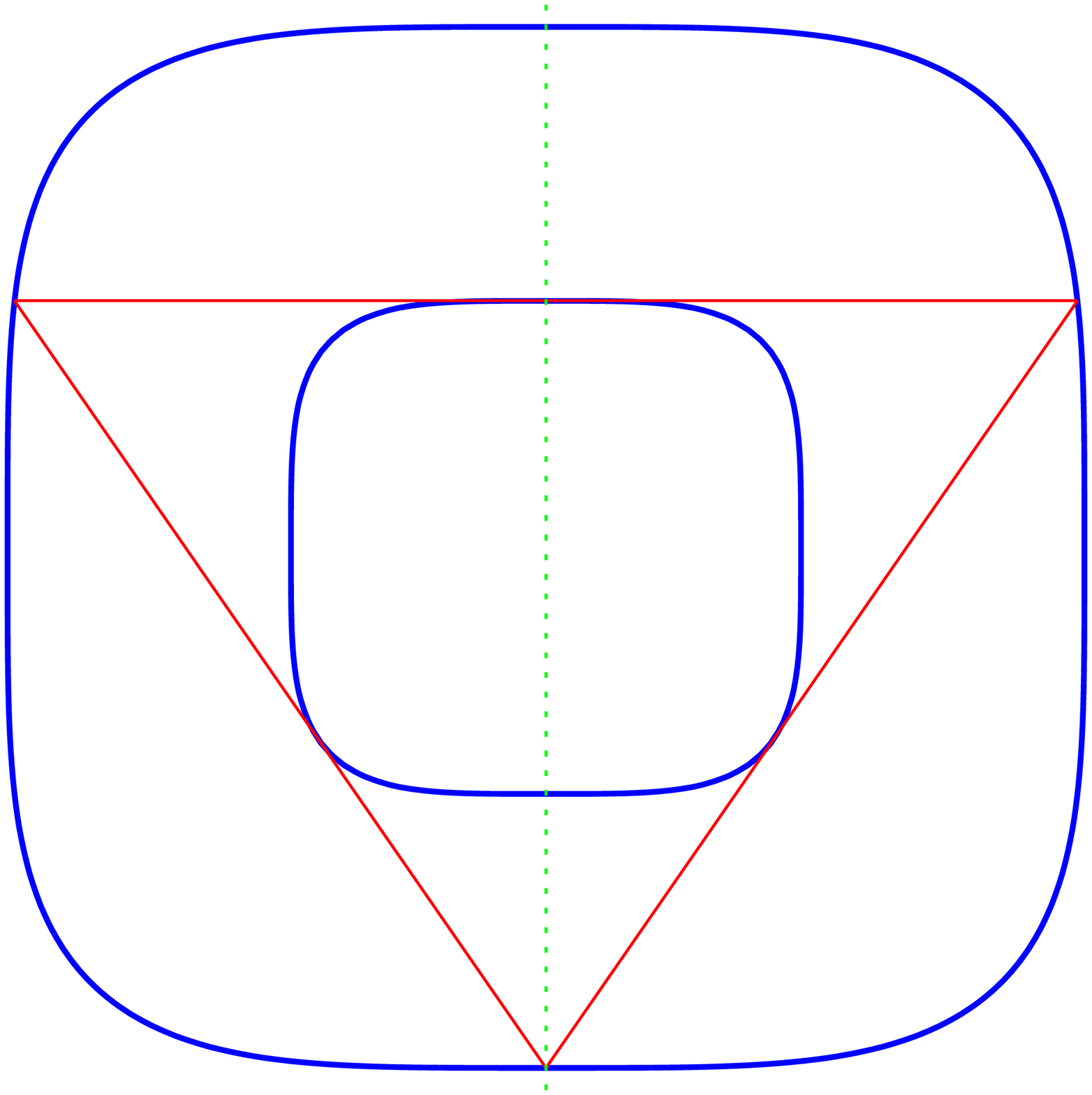,width=0.26\textwidth}$\qquad$
\epsfig{file=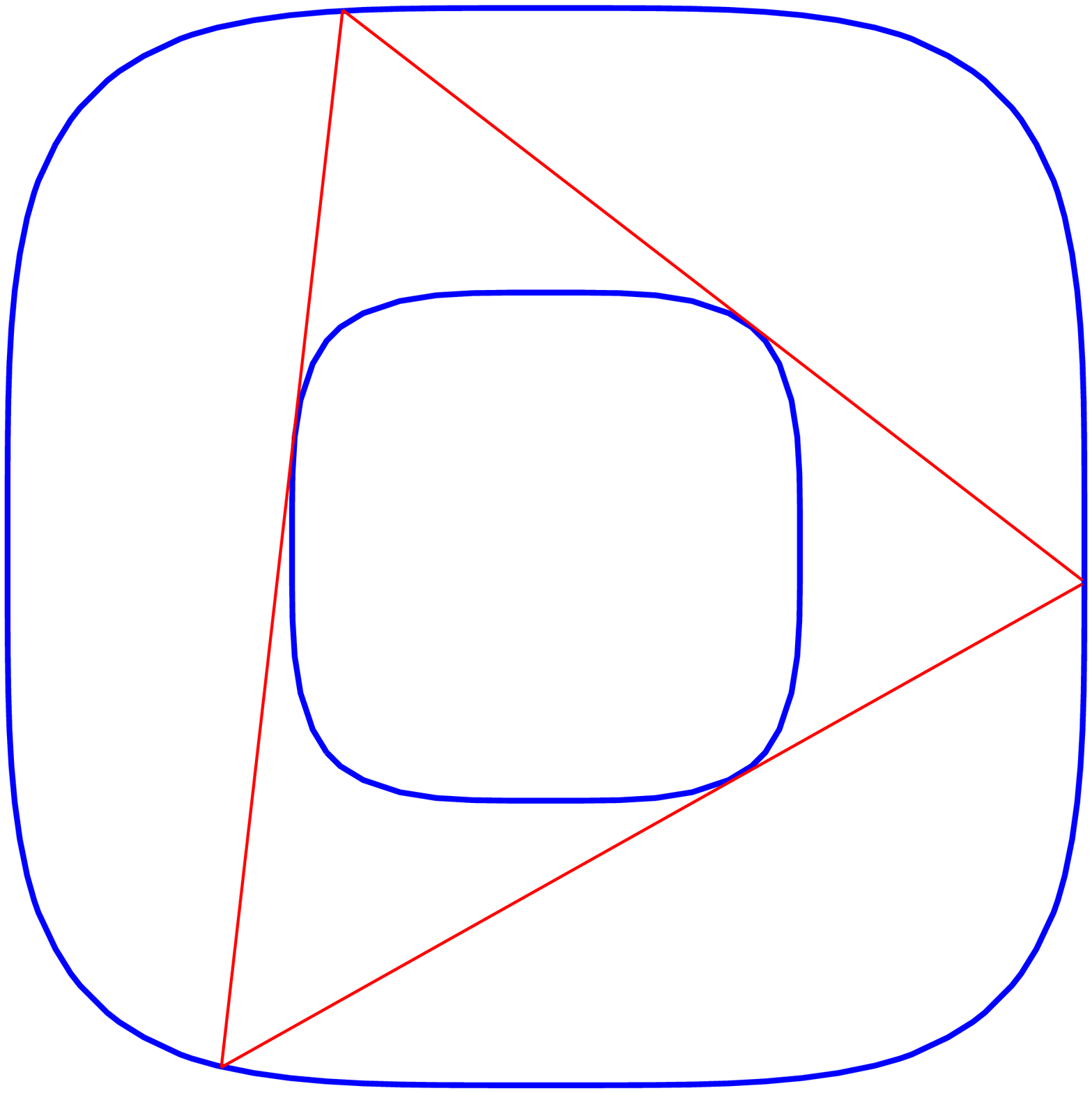,width=0.26\textwidth}$\qquad$
\epsfig{file=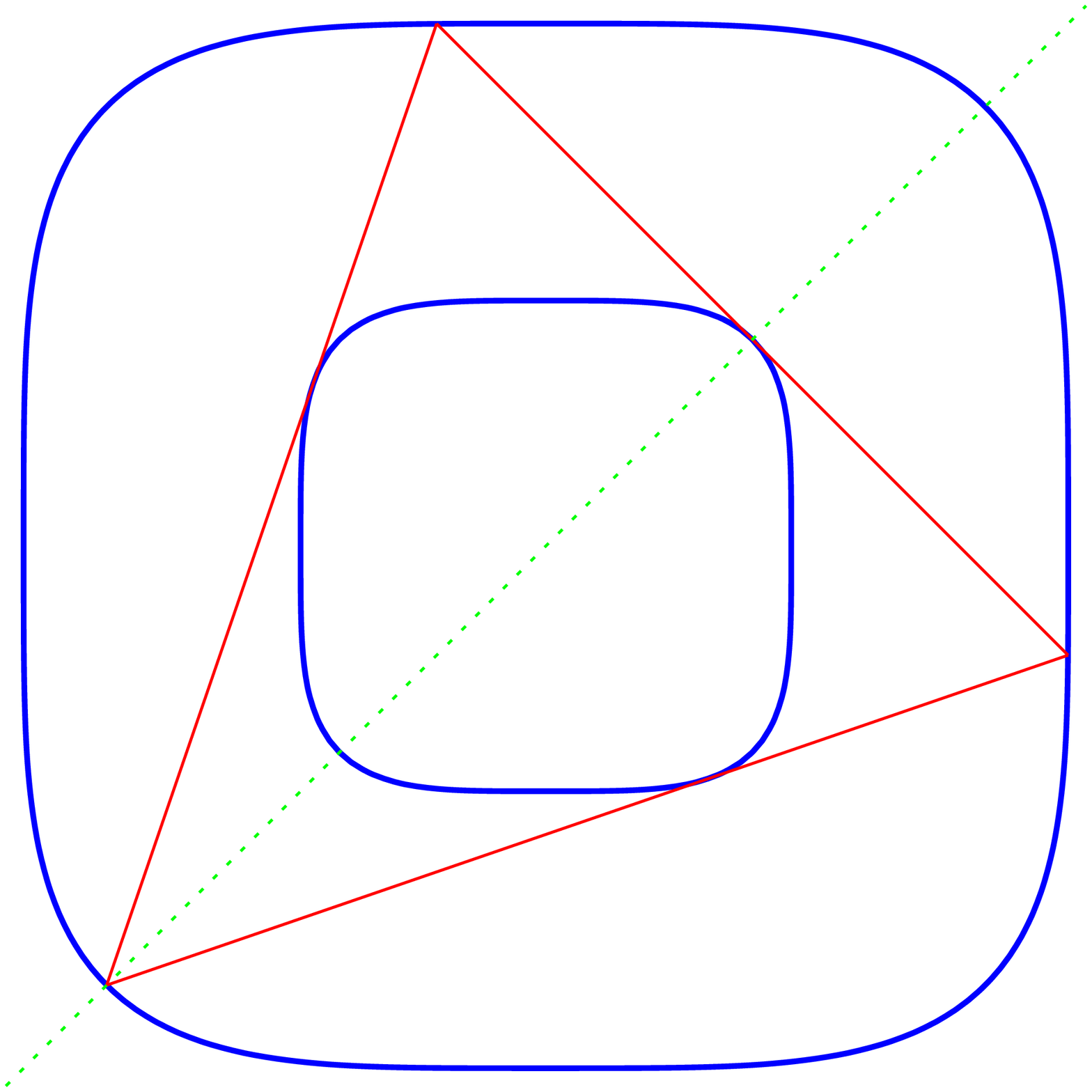,width=0.26\textwidth}
\end{figure}
$\phantom{xx}k\simeq19.8264\phantom{xxxxxxxxxxxxxx}
k\simeq20.1961\phantom{xxxxxxxxxxxxxx} k\simeq20.5087$
\\

\textbf{FIGURE 3.} Three Poncelet's maps with rotation number $1/3$
associated to $\gamma=\{x^4+y^4=1\}$ and $\Gamma=\{x^4+y^4=k\}$.
Notice that the middle one is not symmetric.}
\end{center}

\begin{center}{
\begin{figure}[h]
\hspace{3cm}\epsfig{file=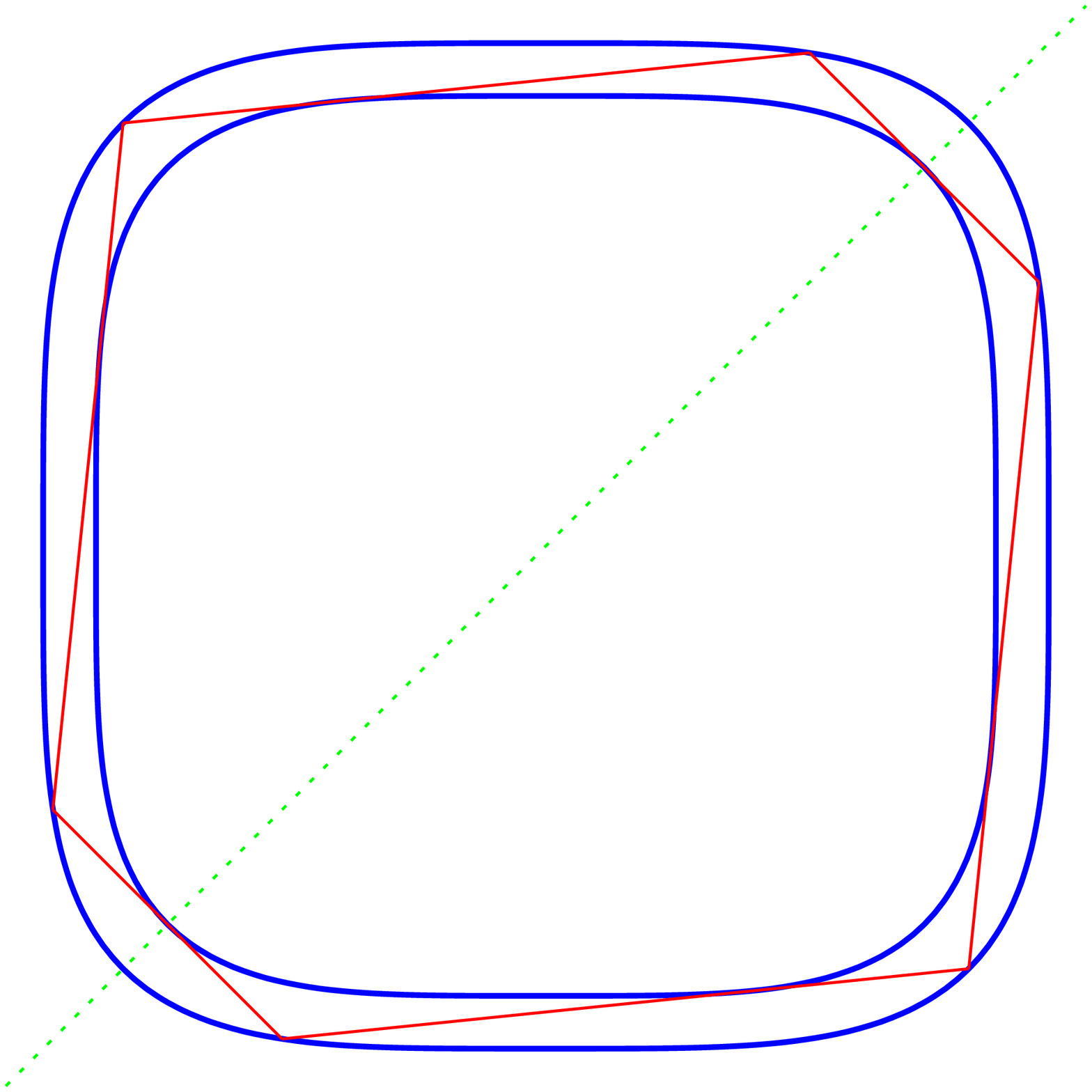,width=0.28\textwidth,}$\qquad\qquad$
\epsfig{file=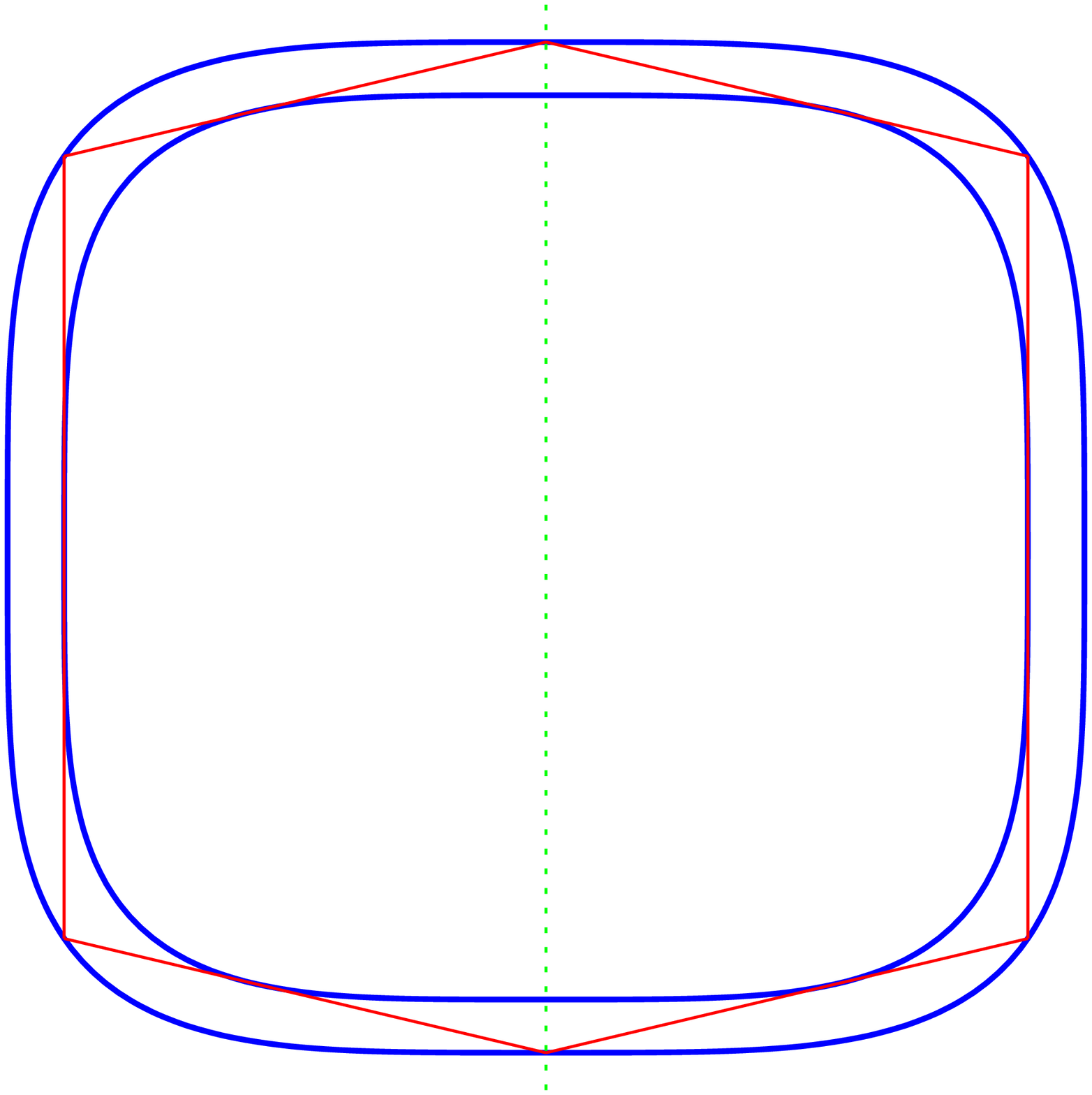,width=0.28\textwidth,} %
\end{figure}
$\phantom{xxxxxx}$ $k \simeq 1.5588$  
$\phantom{xxxxxxxxxxxxxxxxx}$ $k\simeq1.5596$
\\
\vspace{0.4cm}
 \textbf{FIGURE 4.} Two Poncelet's maps with rotation
number $1/6$ associated to $\gamma=\{x^4+y^4=1\}$ and
$\Gamma=\{x^4+y^4=k\}.$}
\end{center}

We have done similar computations to find values of $k$ for which
the corresponding Poncelet's maps have symmetric periodic orbits
(with respect to either the axes or the diagonals) and have rotation
numbers $1/3,$ $1/4$ and $1/6,$ see Figures 3 and 4. All these
results, together with the fact that
\[
\lim_{k\to1} \rho(k)=0  \quad \mbox{and}  \quad \lim_{k\to\infty}
\rho(k)=\frac12,
\]
which simply follows from the geometrical interpretation of $P,$ are
collected in Table 1.

\begin{center}
\begin{table}[h]
\begin{tabular}{lcccccccc}
  \hline
   $k$: &1&$\simeq$ 1.5588& $\simeq$ 1.5596 & 2&8&$\simeq $19.8264&$\simeq$ 20.5087&
   $\infty$\\

   $\rho(k)$&0&1/6 & 1/6 & 1/4 & 1/4 & 1/3 & 1/3&1/2\\
   $\mbox{Symmetry}:$&--&\mbox{diagonals}&\mbox{axes} & \mbox{all} & \mbox{all} & \mbox{axes} & \mbox{diagonals}&--\\
   \hline
\end{tabular}
\caption{Some values of $k$ with rational rotation number and
symmetric period orbits.}
\end{table}
\end{center}

\begin{nota}\label{nota} Notice that, once we have a $n$-periodic orbit
for a given Poncelet's map, then the  six sets obtained by applying
to it either a rotation of $0, \pi/2, \pi$ or $3\pi/2$ radians or
one of the symmetries with respect to the diagonals $\{y=\pm x\},$
are also  $n$-periodic orbits. Moreover, unless the original
periodic orbit has some of the symmetries considered in this
subsection, these six $n$-periodic orbits are all different.
\end{nota}

As a corollary of the above remark we get, for instance, that for
$k=k^*$ the Poncelet map has four 3-periodic orbits. Similarly, for
$k\simeq20.5087$ the corresponding Poncelet's map has four different
periodic orbits, as well. Moreover in both cases it is not difficult
to check that $P$ is not conjugated to the rotation of angle
$2\pi/3.$

\subsection{On the stairs of the rotation function}

Let $\Phi_\lambda$ be  a  smooth one-parameter family of
diffeomorphism of the circle. Let $\rho(\lambda)$ be the rotation
number of $\Phi_\lambda$. Fix a natural number $m.$  Recall that the
existence of an open interval $I_{m}$ where $\rho(\lambda)\equiv
j/m,$ for some natural $j,$ coprime with $m,$   is a  consequence of
the existence of a hyperbolic $m$-periodic orbit for $\Phi_\lambda$
for some $\lambda\in I_{m}.$ For this reason it is said that for
generic one-parameter families of diffeomorphims the graph of the
rotation function has a devil's staircase shape.

The existence of open intervals on which the rotation function
$\rho(k)$ for   Poncelet's maps $P$ associated to $
\gamma=\{x^4+y^4=1\}$ and $\Gamma_k=\{x^4+y^4=k\}$ is constant can
be interpreted  by using the above facts. Here we discuss how the
existence of this type of intervals can also be interpreted more
geometrically.

Consider for instance  the two values of $k$ for which $\rho(k)=1/3$
obtained in  Subsection~\ref{spo},   $k^*\simeq19.8264$ and $\tilde
k\simeq 20.5087.$   By using the method described in
Subsection~\ref{po} we have done a numerical study about the
existence of other values of $k$ having as well a 3-periodic orbit.
We have obtained that for any $k\in [k^*,\tilde k]$ such an orbit
exists, see for instance the middle picture in Figure~3. Moreover
for any $k\in (k^*,\tilde k)$ the orbit that we have found has none
of the symmetries described in this section. By using
Remark~\ref{nota}, we know that for each of these values of $k$, the
Poncelet's map has six different periodic orbits. On the other hand
the boundaries of the interval correspond with values of $k$ for
which some of these six 3-periodic orbits collide giving rise to
some symmetric 3-period orbit, which indeed has to be a multiple
$3$-periodic orbit and so no hyperbolic for~$P.$

We have also checked that a similar phenomenon  occurs when
$k\in[2,8].$ On this interval $\rho(k)\equiv1/4.$

An open problem is to prove if the situation described for $1/3$ and
$1/4$ holds for any rational number in $(0,1/2)$ and also study the
same question for other families of convex ovals.

\subsection{A consequence of the behavior of the rotation
function}\label{ode}

In \cite{CGM} it is proved the following result:

\begin{teo}\label{teor}
Let $\U\subset \mathbb{R}^2$ be an open set and let
$\Phi:\U\rightarrow \U$ be a diffeomorphism such that:
\begin{enumerate}
\item[(a)] It has a smooth  regular first integral $V:\U\rightarrow \R,$ having
its level sets $\Gamma_k=:\{p\in\U\,:\, V(p)=k\}$ as simple closed
curves.
\item[(b)] There exists a smooth function $\mu:\U\rightarrow \R^+$ such
that for any $p\in \U,$
\begin{equation}\label{mu}
\mu(F(p))=\det(D\Phi(p))\,\mu(p).
\end{equation}
\end{enumerate}
Then the map $\Phi$ restricted to each $\Gamma_k$ is conjugated to a
rotation with rotation number $\tau(k)/T(k)$, where $T(k)$ is the
period of $\Gamma_k$ as a periodic orbit of the planar differential
equation
\[
\dot p=\mu(p)\left(-\frac{\partial V(p)}{\partial
p_2},\frac{\partial V(p)}{\partial p_1}\right)
\]
and $\tau(k)$ is the time needed by the flow of this equation for
going from any $q\in\Gamma_k$ to $\Phi(q)\in\Gamma_k.$
\end{teo}

Notice that it provides a way to check whether integrable planar
maps $\Phi$ of the circle are conjugated to rotations or not. It
consists in studying the existence of solutions $\mu$ of the
functional equation \eqref{mu}.

A natural problem in this context is to study under which conditions
over an integrable map $F,$ the functional equation \eqref{mu} has
non-trivial solutions. Let us see  that  the Poncelet's map $\Phi=P$
constructed in Section~\ref{sec1} associated to the curves
$\{x^4+y^4=1\}$ and $\{x^4+y^4=k\}, k>1,$ and defined in the set
$\mathcal{U}=\{2<x^4+y^4<8\},$ provides an example of map $F$ for
which equation \eqref{mu} has no solution. Notice that $P$ is
clearly integrable, with first integral $V(x,y)=x^4+y^4.$ If
associated to $P$ it would exist a function $\mu$ satisfying the
functional equation \eqref{mu} then, by the results of the previous
subsection and Theorem~\ref{teor},   $P^4=\mbox{Id}$ on $\U,$ result
that is trivially false.

\vspace{2cm}

\section{A new proof of Poncelet's Theorem}\label{newproof}

This section is devoted to give a new proof of Poncelet's Theorem
based on Theorem~\ref{teor}. To do this it is more convenient to
take coordinates in such a way that the outer ellipse is given by
the circle $\G=\{x^2+y^2=1\}$ and the inner one is given by the set
$\g=\{g(x,y):=A{x}^{2}+  Bxy+ C{y}^{2}+Dx+Ey+F=0\}.$ Moreover it is
not restrictive to assume that $g$ is positive on the exterior of
the ellipse. Following our point of view we will consider the
Poncelet's map $P$ defined from the region $\mathcal{U}=\{(x,y)\in
\mathbb{R}^2\,:\, |(x,y)|>d_0\}$ into itself, where $d_0$ is the
maximum distance between $\g$ and the origin. Hence
$\G\subset\mathcal{U}$. Some tedious but straightforward
computations, done with an algebraic manipulator, give that the
Poncelet's map writes as:

\begin{equation}\label{pon2d}
P( {x,y} )=\left({\frac {-N_{{1}}N_{{2}}-4 N_{{3}}\sqrt
{\Delta}}{M}}, {\frac {-N_{{1}}N_{{3}}+4 N_{{2}}\sqrt
{\Delta}}{M}}\right)
\end{equation}
where

\noindent $\begin{array}{rl}
    N_1= & 4 AF+4 CF-{D}^{2}-{E}^{2}+ 2\left(2 C{D} - BE
 \right) x+2\left( 2 AE- {D} B \right) y\\&+\left(4 AC -{B}^{2}\right)\left[ x^2+
{y}^{2}\right],
  \end{array}
$

\noindent$
\begin{array}{rl}
    N_2= & \left( 4 CF+{{D}}^{2}-4 AF-{E}^{2} \right) x+ 2\left(
 DE-2 BF \right) y\\&
 + \left[2\left(2 C{D} - BE \right)
+ \left( 4 AC-{B}^{2}
 \right)x\right]\left[{x}^{2}+{y}^{2}\right],
  \end{array}
$

\noindent$
\begin{array}{rl}
   N_3= & 2\left(  DE-2 BF \right) x+ \left( {E}^{2}-{{
D}}^{2}-4 CF+4 AF \right) y\\&+ \left[2\left( 2 AE- B{D} \right)
 +\left( 4 AC-{B}^{2}
 \right) y\right]\left[{x}^{2}+ {y}^{2}\right],
  \end{array}
$

\noindent$
\begin{array}{rl}
M=&{E}^{4}+{{D}}^{4}+16 {F}^{2}{A}^{2}-16 EFB{D}+16 {F}^{2}{
C}^{2}+16 {F}^{2}{B}^{2}+2 {{\it
D}}^{2}{E}^{2}-8 FA{{D}}^{2}\\
& +8 FA{E}^{2}-8 F{E}^{2}C+ 8 FC{{D}}^{2}-32 {F}^{2}AC +
\left( - 12 {{D}}^{2}EB-16 AEFB\right.\\
&\left. -32 FCA{D}+4 B{E}^{3}+ 16 {D}  A{E}^{2}+8 C{{\it D}}^{3}-8
C{\it
D} {E}^{2}-16 FEBC\right.\\
&\left. +32 F{C}^ {2}{D}+16 {D} {B}^{2}F \right) x+\left( -32 FACE-8
A{{  D}}^{2}E+  16 C{{D}}^{2}E+
16 {B}^{2}EF\right.\\
&\left. +32 F{A}^{2}E+4 B{{  D}}^{3}+8 {E}^{3}A-16 FA{\it D} B-16
BC{D} F-12 B{ \it
D} {E}^{2} \right) y\\
&  +\left( 16 {C}^{2}{{D}}^{2}+6 {E}^{2
}{B}^{2}+8 FA{B}^{2}-8 F{B}^{2}C-32 F{A}^{2}C+32 F{C}^{2}A\right.\\
&\left. +2 {B}^ {2}{{D}}^{2}+16 {A}^{2}{E}^{2}-16 {\it D} BEC-8
A{E}^{2}C+8  AC{{D}}^{2}-16 AB{D} E \right)
{x}^{2}\\
& +\left( 32 AC{  D} E-8 {D} {B}^{2}E-64 BCAF+16 {B}^{3}F \right)
yx+
 \left( -16 {D} BEC+16 {A}^{2}{E}^{2}\right.\\
&+\left. 2 {E}^{2}{B}^{2}-8 FA{ B}^{2}-32 F{C}^{2}A+6 {B}^{2}{{
D}}^{2}+16 {C}^{2}{{D}}^{2 }+8 A{E}^{2}C-16 AB{\it
D} E\right.\\
&\left. +8 F{B}^{2}C+32 F{A}^{2}C-8 AC{{  D}}^{2} \right)
{y}^{2}\\
& +4(4AC-B^2)\left[(2CD-BE)x+(2AE-BD)y\right]\left[x^2+y^2\right]
\\
&+  \left( 4AC-B^2 \right)^2 \left[{x}^{2}+y^2\right]^2,
  \end{array}
$

\noindent and

\noindent$\begin{array}{rl} \Delta=&\left(A{E}^{2}-BDE
+C{{D}}^{2}+F({B}^{2}-4 AC)
 \right) \left[
A{x}^{2}+  Bxy+ C{y}^{2}+Dx+Ey+F\right].
 \end{array}$

\vspace{0.2cm}

In Figure  5 we show some points of two orbits generated by this map
corresponding to different ellipses. Recall that by construction the
map $P$ given in \eqref{pon2d} is a diffeomorphism on $\mathcal U$
and $V(x,y)=x^2+y^2$ is a first integral for it.

Then in order to apply Theorem~\ref{teor} we only need to find a
smooth function, defined on $\mathcal U,$ such that the functional
equation \eqref{mu} with $\Phi=P$ holds. It can be seen, again by
using an algebraic manipulator,  that a function satisfying this
equality
 is given by
\[
\mu(x,y)=\sqrt{(x^2+y^2)g(x,y)}=\sqrt{(x^2+y^2)\left(A{x}^{2}+  Bxy+
C{y}^{2}+Dx+Ey+F\right)}.
\]
Hence Poncelet's Theorem follows. Notice that our proof also works
when the rotation number of the Poncelet's map is  irrational.

We end this section with some comments of how we have got the above
function $\mu.$ By using the Change of Variables Theorem it is not
difficult to check that the existence of a positive function $\mu$
satisfying equality \eqref{mu} implies that the absolute continuous
measure
\[
\nu(\mathcal{B}):=\iint_{\mathcal B} \frac 1{\mu(x,y)}\,dxdy,
\]
is an invariant measure for $P$, that is
$\nu(P^{-1}(\mathcal{B}))=\nu(\mathcal{B})$ for any measurable set
$\mathcal{B}\subset\mathcal{U}.$ On the other hand, one of the
proofs given in \cite{T} --the one of \cite{K}-- geometrically
constructs an invariant measure on the outer circle to prove
Poncelet's Theorem. Inspired on this construction we have been able
to extend this mesure to the whole $\mathcal U$ and as a consequence
we have got a suitable $\mu.$

\vfill

\begin{center}{
\begin{figure}[h]\hspace{2cm}
\epsfig{file=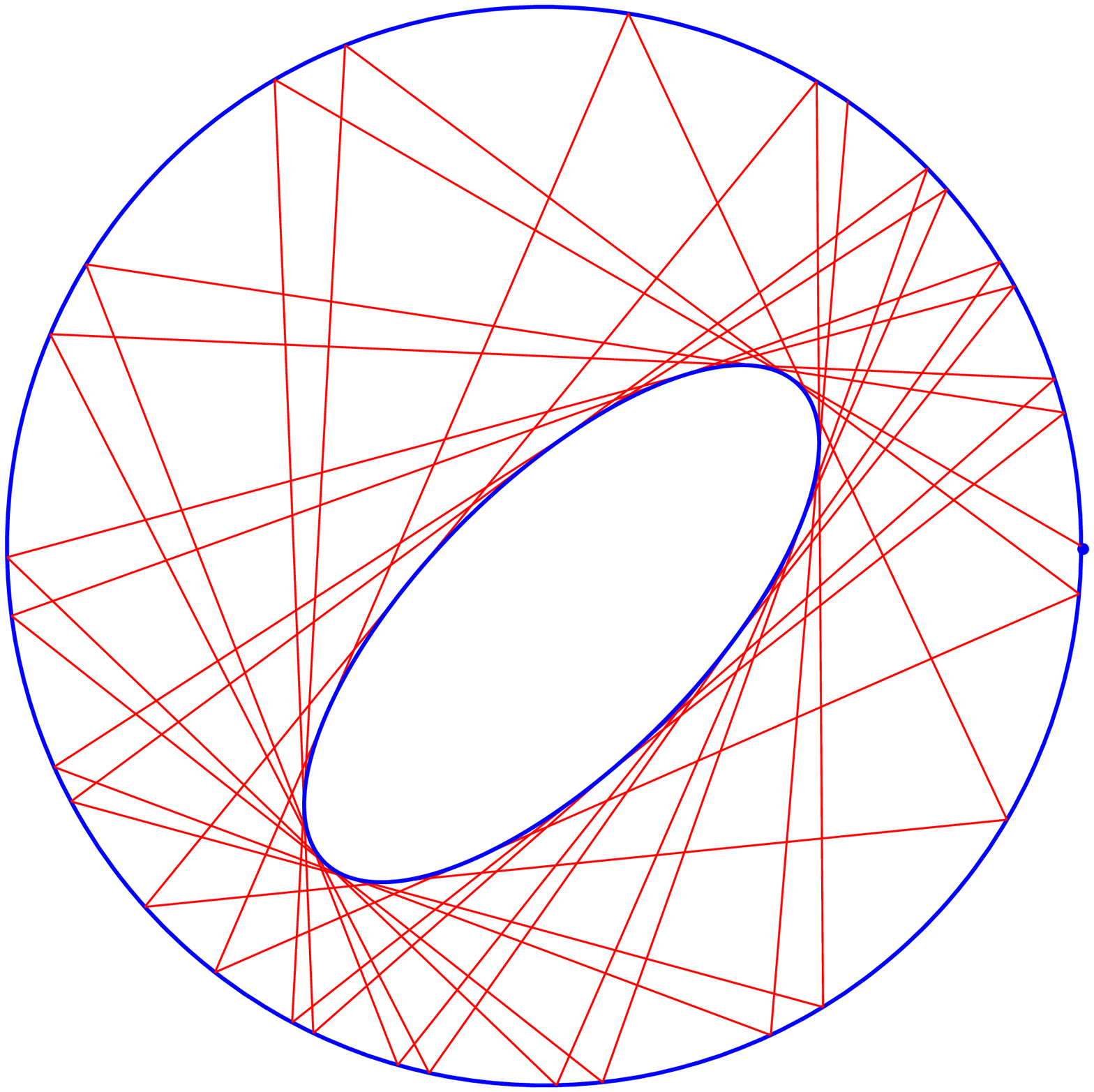,width=0.35\textwidth,}$\qquad\qquad$
\epsfig{file=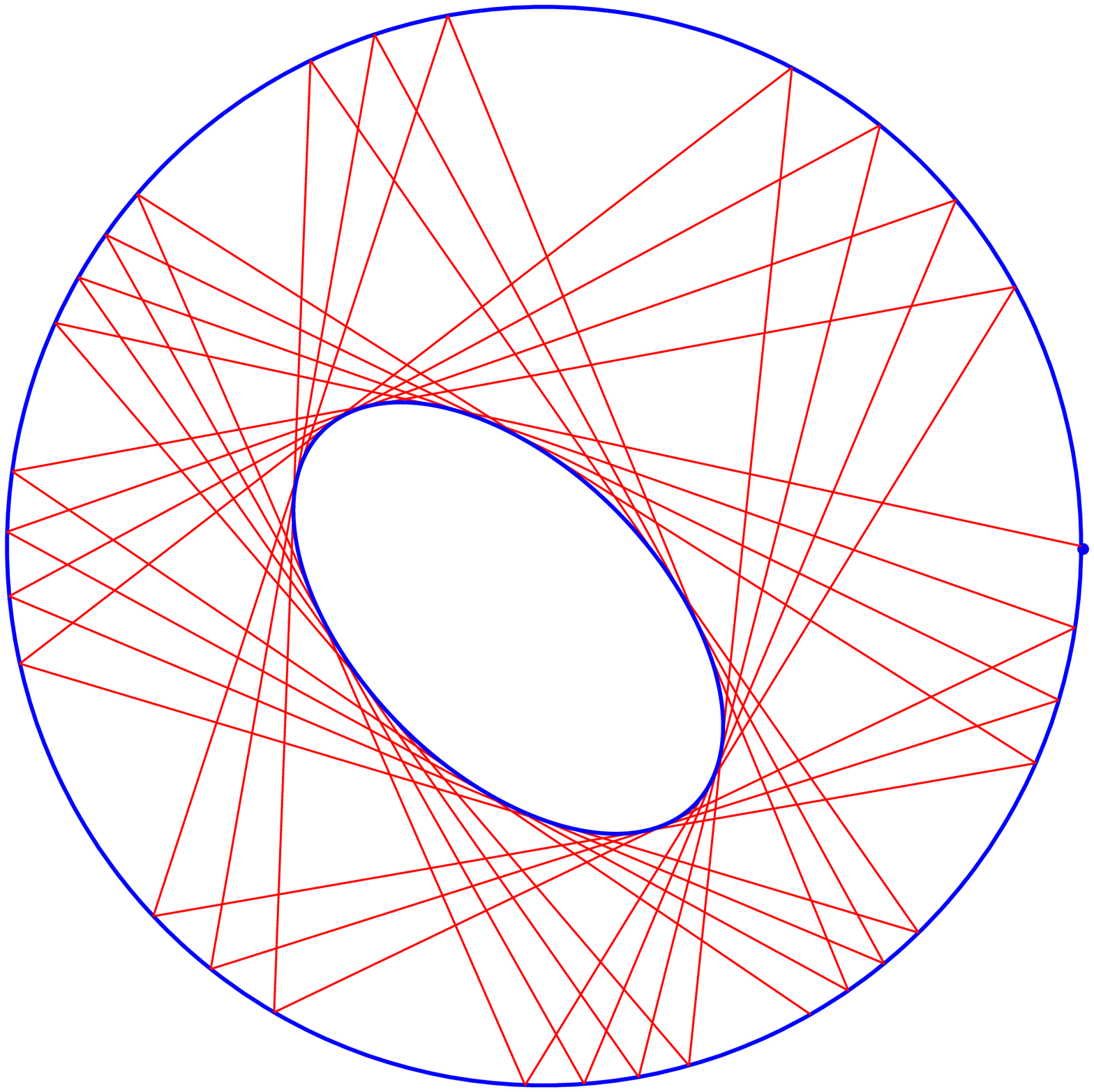,width=0.35\textwidth,}
\end{figure}
\textbf{FIGURE 5.} Thirty points of two Poncelet's orbits. \\The dot
corresponds to the initial condition.}
\end{center}

%
%
%
%
%

%
%
%
%

\noindent {\bf Acknowledgments.} The authors are  supported by
DGICYT through grants MTM2005-06098-C02-01 (first and second
authors) and DPI2005-08-668-C03-1 (third author). They are also
supported by the Government of Catalonia through some SGR programs.

The second author thanks Emmanuel Lesigne for stimulating
discussions about  Poncelet's maps.

\end{document}